\documentclass[letterpaper, 11pt]{article}
\usepackage{amsmath, amssymb}
\usepackage{bm}
\usepackage{graphicx}
\usepackage[font={small,it}, labelfont=bf]{caption}                
\usepackage{subcaption}
\usepackage[margin=1in]{geometry}                         
\usepackage{setspace}    
\usepackage{wrapfig}                                
\usepackage{xcolor}
\usepackage{multirow,tabularx}
\newcolumntype{Y}{>{\centering\arraybackslash}X}

\usepackage{lineno}

\usepackage{cite}

\usepackage{titlesec} \titleformat*{\section}{\bfseries}              
\titleformat*{\subsection}{\bfseries}                       

\usepackage{titling}
\pretitle{\begin{flushleft}\bfseries\vskip -5em\Large}
\posttitle{\end{flushleft}}
\preauthor{\begin{flushleft}\vskip -0.5em \small}
\postauthor{\end{flushleft}}
\predate{}
\postdate{\vskip -1.5em}

\usepackage{authblk}

\author[1,*]{Kristofer G. Reyes}
\author[2]{Jiaqian Liu}
\author[3]{Carlos Juan D\'{i}az Vargas}

\affil[1]{Department of Materials Design and Innovation, University at Buffalo.
Buffalo, NY, 14260}
\affil[2]{Department of Computer Science, University at Buffalo.  Buffalo, NY, 
14260}
\affil[3]{Department of Computer Science, University of Puerto Rico. San Juan, 
PR, 00925}

\affil[*]{Corresponding Author: kreyes3@buffalo.edu}

\title{Decision-Making Under Uncertainty for Multi-stage Pipelines:
Simulation Studies to Benchmark Screening Strategies}
\date{}

\begin{document}
\maketitle

\begin{abstract}
Multi-stage screening pipelines are ubiquitous throughout experimental and
computational science. Much of the effort in developing screening pipelines
focuses on improving generative methods or surrogate models in an attempt to
make each screening step effective for a specific application. Little focus has
been placed on characterizing a generic screening pipeline's performance with
respect to the problem or problem parameters. Here, we develop methods to
codify and simulate features and properties about the screening procedure in
general. We outline and model common problem settings, and identify potential
opportunities to perform decision-making under uncertainty for optimizing the
execution of screening pipelines. We then illustrate the developed methods
through several simulation studies. We finally show how such studies can
provide a quantification of the screening pipeline performance with respect to
problem parameters, specifically identifying the significance of stage-wise
covariance structure. We show how such structure can lead to qualitatively
different screening behaviors, and how screening can even perform worse than
random in some cases.
\end{abstract}

\section{Introduction}

Multi-stage screening pipelines are ubiquitous throughout experimental and
computational science. In such a pipeline, a set of initial candidates are
passed through several evaluative stages. Based on the collective performance
of the candidates at any particular stage of the pipeline, a subset of
candidates are selected to advance to the next stage; the rest are ``screened
out". Typically, this screening is done to minimize how many candidates are
evaluated throughout the pipeline, with the goal of identifying promising or
optimal candidates with respect to some experimental objective. Often, later
stages are more costly, in time, money or other measures of effort.

For example, experimental and virtual screening are used in drug-design
\cite{andricopulo2008virtual, SCHNEIDER200264, kandeel2020virtual,
bajorath2002integration}, chemistry \cite{shoichet2004virtual}, material
science\cite{gomez2016design, pyzer2015high, pyzer2016bayesian}. Here,
candidate drugs, molecules or materials designs are first tested using
preliminary computer models, and promising candidates are then tested in
real-world experiments. Another example is multi-modal characterization of
materials. In this context, a material sample is characterized using one of
several  techniques, from atomic force microscopy available in many labs, to
X-ray scattering at a synchrotron. Here, a materials scientist may employ
an \emph{ad hoc} screening procedure to identify which materials sample candidates
are passed through this characterization pipeline. As autonomous materials
platforms begin to incorporate more complex experimental or hybrid
computational/experimental structure \cite{abdel2021self, gongora2021using,
baek2021problem}, understanding how screening pipelines could operate in an
autonomous manner could allow such platforms the ability to make optimal
screening pipeline decisions as well.

The use of a screening procedure relies on a heuristic belief that the screening
pipeline will effectively narrow down a search space to a small set of promising
leads. But there is no generic guarantee this will work. Consider, for example
the case of identifying promising drugs, where the ``gold-standard" measure for
effectiveness occurs at the last stage of a screening pipeline, which could be
performing physical trials of the drug. Imagine the first stage is a machine
learning (ML) model meant to predict this effectiveness based on a candidate
drug's descriptors, but its predictions are only slightly correlated or even
anti-correlated with the gold-standard measure. In this case, the first stage
of the pipeline may aggressively and prematurely screen out good candidates. 

To mitigate against this threat, much of the development around screening
pipelines has been the development of stages and models that are more strongly
correlated to the experimental objectives, either through better surrogate
models\cite{tropsha2007predictive, azencott2007one, ain2015machine,
ballester2010machine} or through better generative
models\cite{sanchez2018inverse, lim2018molecular, gao2020synthesizability} that
determine the initial candidate pool to begin with. Yet fundamental and generic questions
persist: What aspects of a problem make the screening procedure effective?
Can we quantify the effectiveness of the screening procedure?  How do we deal
with the uncertainties behind the relationships between stages? Are our current
pipeline strategies or \emph{ad hoc} \emph{decision-making policies} sufficient
in executing an effective pipeline?  Coming back to the drug-design example
above, we may wish to quantify how poorly the screening procedure performs as a
function of the anti-correlation between the first and final stages.

To answer such questions, we must develop a theory where the general pipeline
itself is the subject of study, specifically from a perspective of uncertainty
quantification and decision-making under uncertainty.  In this paper we outline
a formalization of a family of screening pipeline problems, and
simulation-based policies that we use to make decisions before and during the
screening procedure. To accomplish this, we build a model of statistical
structure between candidates and stages, and offer algorithms for efficient
simulation-based decision-making based on the assumed statistical
structure model. Lastly, we present some simulation results to illustrate how
such models, simulations, and policies can be used to characterize and optimize
screening pipelines.

\section{Problem Definitions}

We consider a screening procedure with $m$ initial sampled candidates, $x_1,
x_2, ..., x_m$. These candidates can represent drugs, molecules, materials, or
even synthesis recipes. Often, these initial candidates are obtained via
generative models, ensemble methods, or evolutionary techniques -- methods
capable of proposing a batch of candidates to test.

We will also consider $n$ stages $f_1, f_2, .., f_n$, which we shall view as
functions mapping candidates $x_i$ to corresponding scalar quantities of
interest $Y_{i,j} = f_j(x_i)$ as measured by the $j$-th stage. For example in
drug design, the stages could represent molecular property predictors, docking
score calculations, density functional theory or molecular dynamics
simulations, or physical trials.  The outcomes for these stages could be
computed binding efficacy, free energy calculations, or measures toxicity and
synthesizability.  In manufacturing, the stages could include low and high
fidelity finite element simulations, rapid physical prototyping via 3D
printers, and fabrication of a final deliverable.  The outcomes for these
stages could be measures of mechanical toughness or resiliency.

In the abstract, we call these values \emph{scores}, so that $Y_{i,j}$ is the
score assigned to candidate $x_i$ at stage $j$.  We  interpret these scores as
the averaged ground truth values. That is, if we were to repeatedly evaluate a
candidate $x_i$ through the $j$-th stage, the average of the set of obtained
scores will approach the number $Y_{i,j}$.  Thus the $Y_{i,j}$ represents the
score of candidate $x_i$ through the $j$-th stage, averaging over noise or any
inherent randomness implicit in making the appropriate measurement and
obtaining the quantity of interest. 

Throughout, we assume we do not know the scores prior to evaluating the
candidates through each stage.  As \emph{a priori} uncertain quantities, we
treat the scores within a \emph{Bayesian} context by modeling such scores using
a multivariate probability distribution. That is, throughout the screening
procedure, we can maintain a probability distribution for the scores $Y_{i,j}$
that captures predicted estimates of scores not yet evaluated as well as an
uncertainty quantification of those predictions and an estimated correlation
between different scores. As more data is obtained through the screening
procedure, we can update this probability distribution through Bayes's Law.

\subsection{Objectives}

These scores allow us to define two objectives. First is \textbf{stage-wise
multi-fidelity optimization}, in which we view each stage as measuring the same
quantity, but with variable accuracies (i.e.  fidelities) and costs.  In
general, we will assume that the last stage $f_n$ is the ground-truth or
gold-standard that we ultimately care about and is the most expensive. The
prior stages exist to perform cheaper screening of candidates. For example,
$f_n$ could be performing a real experiment, while prior stages could be a
sequence of increasingly accurate (and increasingly costly) simulations. A
common pipeline in this setting consists of three stages: $f_1$ is an ML
model, $f_2$ is a physics-based simulation, while $f_3$ is a physical
experiment.  Each stage could produce as an output some figure of merit, such
as some material property to optimize, and the goal is to find the candidate
that optimizes this property in the real-life experiment. Another similar
setting is when stages correspond to different materials characterization
methods, where the final such method could be an costly trip to a beamline. The
stage-wise multi-fidelity optimization problem shall be the main focus of this
paper.

Another objective is \textbf{stage-wise multi-objective optimization.} In this
setting, we view each stage as characterizing different quantities of interest.
For any given candidate $x_i$, the scores $Y_{i,1}, Y_{i,2}, ..., Y_{i, n}$
encode a set of different characterizations of candidate $x_i$'s material
properties or structure, for example. Under the multi-objective objective, we
want to identify those candidates $x$ that either 1) optimize some scalar
utility $u(x) = u(f_1(x), f_2(x), ..., f_n(x))$, or 2) lie on the
$n$-dimensional Pareto frontier\cite{deb2014multi}.

Both objectives can be characterized by defining the \emph{final reward}
obtained at the end of the screening procedure. For any screening procedure we
obtain a subset $\mathcal X_S \subseteq \left\{ x_i \right\}_{i=1}^m$ of
candidates that ``survived" the screening process. For the multi-fidelity
objective, we may define the reward as the maximum of the final scores for all
such screened candidates, $R = \max_{\mathcal X_S} f_m(x)$. For
multi-objective optimization, the final score could be either the maximum
scalar utility over all surviving candidates, $R = \max_{\mathcal X_S}
u(x)$, or some measure of improvement to the Pareto front, such as a
hypervolume improvement\cite{BEUME20071653}, $R = H(\mathcal X_S) - H_0.$ Here
$H(\mathcal X_S)$ and $H_0$ are respectively the hypervolumes defined by the
Pareto front generated by the surviving points $\mathcal X_s$ and prior data
points, and the Pareto front obtained from just the prior data points.

\subsection{Pipeline policies for the optimal allocation design problem}
Regardless of how we define the reward, the goal in executing the pipeline
is to do so in a way that maximizes the realized reward on average. Throughout,
we shall refer to how we execute the pipeline as a \emph{pipeline policy}. The
specific decisions we make during the execution of the pipeline will ultimately
impact the rewards we gain at the end of the pipeline.  Examples of such 
decisions include: 1) how many or what proportion of candidates may pass from
one stage to the next, 2) how costly each stage is (to the extent that this
cost is possible to tune) and 3) how many or  which initial candidates do we
generate. 

In this paper, we shall decompose a pipeline policy into two policies. The 
\emph{inter-stage policy} selects which candidates pass to the next stage 
based on the scores they receive at the current stage. In contrast, the 
pipeline \emph{meta-policy} select hyperparameters and constraints that 
the inter-stage policy must satisfy. The main example used throughout this paper 
is selecting an \emph{allocation} $\mathbf m = (m_1, ..., m_n)$ via a 
meta-policy.  An allocation specifies how many candidates are evaluated at 
each stage.  In this notation, $m_1 = m$ candidates are evaluated in the 
first stage, then $m_2$ candidates are selected (by the inter-stage policy) for
the second stage, and so on. Throughout, we require each stage to screen
out at least one candidate, and that some candidates survive to the end.
That is, we require $m_1 = m > m_2 > \cdots > m_n \geq 1.$

The meta-policy selects allocations, while the inter-stage policy selects
candidates given an allocation. This distinction is not fundamental. We can
imagine an inter-stage policy that selects both the number of candidates to
pass trough to the next stage, and which specific candidates do so. However for
this paper, we focus on the decomposition of decision variables described
above. We have previously studied a similar nested-batch decision-making policy
in the case of Bayesian optimization for the co-design of materials and
devices\cite{wang2015nested} and have found the decomposition computationally
favorable in practice.

We employ this nested decision policy to solve the \textbf{optimal allocation
problem}. Here, we must select the allocation $\mathbf m^\star$ that is likely
to result in maximal rewards under a specific inter-stage policy. The
meta-policy performs this selection prior to the actual execution of the
pipeline. Instead, it relies  on limited prior knowledge about the scores
to design an allocation. We encode this prior knowledge using a
probability distribution, a \emph{(Bayesian) prior} on the scores. This prior
captures uncertainties about the scores and the assumed statistical
relationship (the \emph{covariance}) between them.

To design an allocation using a prior, we will \emph{simulate} the pipeline
under some inter-stage policy. That is, we can sample a ground truth $\mathbf y
= (y_{i,j})$ from the prior, and use this sample to execute the screening
pipeline. That is, scores are obtained from the sampled ground truth, $f_j(x_i)
= y_{i,j}$, perhaps with some noise added in. Using these scores, the
inter-policy screens out candidates for the next stage, and so on.  Executing
the pipeline this way yields a sample of rewards for the given sampled score
and allocation, $r(\mathbf y, \mathbf m)$.  Repeating the simulations under
different samples of ground truth yields a representative sample of rewards. We
treat this sample as an empirical estimate of the distribution of the reward
$R(\mathbf m)$ for a given allocation $\mathbf m$.  Given this distribution for
each allocation $\mathbf m$, the meta-policy selects an allocation $\mathbf
m^\star$.

\subsection{Meta-policies for closed-loop versus single-shot settings} The
distribution on rewards $R(\mathbf m)$ is a quantification of uncertainty. The
meta-policy may consider this uncertainty, especially  when in a
\textbf{closed-loop setting}. In this setting, the meta-policy selects an
allocation $\mathbf m^\star_1$ based the reward distribution that the prior
distribution induces.  We then implement this decision, running the pipeline in
real-life subject to this allocation. Doing so yields observations of the true
scores, which we use to update the Bayesian prior. The updated prior induces a
new reward distribution, and so the meta-policy can select a second allocation
$\mathbf m^\star_2$.  This procedure is iterated until some total budget is
expended. If the meta-policy is effective, by the end of the campaign, it will
have identified a near-optimal allocation.

In this closed-loop setting, we can use relevant policies from
Bayesian Optimization \cite{snoek2012practical}. They acknowledge the 
uncertainties present in the reward distribution $R(\mathbf m)$ by
balancing between exploration (suggesting random $\mathbf m$ or the one
corresponding to the distribution $R(\mathbf m)$ with the most variance) and
exploitation (suggesting the $\mathbf m$ that maximizes the average, or
expected, value of $R(\mathbf m)$).  For example, the Upper Confidence Bound
(UCB) policy \cite{auer2002finite} explicitly encodes this balance by selecting
the allocation $\mathbf m$ that maximizes the value:
\[A_{\text{UCB}}(\mathbf m) = \mathbb E\left[R(\mathbf m)\right] + c \cdot \sqrt{\text{Var}\left[R(\mathbf m)\right]},\]
where $\mathbb E[R(\mathbf m)]$ is the expected value of the reward 
distribution and $\text{Var}\left[R(\mathbf m)\right]$ is the variance of the
distribution. The square-root of the variance is the standard deviation, which
measures the spread of the distribution.  The expectation term favors
allocations that have a high expected value (exploitation), while the variance
term favors allocations with more reward uncertainty (exploration). The
parameter $c$ balances the two terms.

Such policies are truly effective if we can execute several iterations of the 
closed-loop. They are limited in the case where a budget reduces the number of 
times the entire pipeline can be executed to a few iterations. In the extreme
\textbf{one-shot setting}, the pipeline can be executed once. In this case,
we do not have the luxury to be explorative in our selection of $\mathbf m$. 
Hence the only reasonable meta-policy in the one-shot setting is pure 
exploitation (XPLT):
\[A_\text{XPLT}(\mathbf m) = \mathbb E\left[R(\mathbf m)\right].\]
While this limits our choice of the meta-policy, we may still consider
different inter-stage policies to guide the single execution of the pipeline.
However, for this paper, we shall focus explicitly on the one-shot setting and
the XPLT meta-policy. Yet we note that the development and analysis of proper
inter-stage policies specifically for multi-stage pipelines could be a fruitful
topic of study.

\subsection{The cost-constrained problem}

In practice, later stages of the pipeline can be significantly more expensive
than earlier ones. They represent real experiments, high-fidelity simulations,
or material characterization at a high-demand user facility.  Let $\mathbf c =
(c_1, c_2, \hdots, c_n)$ be a cost vector so that $c_j$ is the cost to evaluate
a candidate at stage $j$.  We assume the costs are additive, and that we cannot
amortize costs or parallelize evaluations.  Therefore, the costs to evaluate
$m_j$ candidates in stage $j$ is simply $m_jc_j$, and the total cost $C(\mathbf
m)$ for an allocation is $C(\mathbf m) = \mathbf m \cdot  \mathbf c$, which we
require to be less than some budget $\mathbf C_\text{max}$.

We frame the XPLT meta-policy via the (linearly) constrained integer program: 
\[
\begin{aligned}
& \underset{\mathbf m \in \mathbb Z^n}{\text{maximize }} & & A_\text{XPLT}(\mathbf m)  \\
& \text{subject to: } & & m_1 = m > m_2 > \cdots m_n \geq 1 \\
& & & \mathbf m \cdot \mathbf c \leq C_\text{max}
\end{aligned}
\]
The meta-policy selects the allocation that solves this program.  While 
the function $A_\text{XPLT}$ is non-linear in $\mathbf m$, it is monotone with 
respect to the dominance partial-order on points  in $\mathbb Z^n$. That is, 
for two allocations $\mathbf m = (m_1, \hdots, m_n)$ and 
$\mathbf m^\prime = (m^\prime_1, \hdots, m^\prime_n)$ such that 
$m^\prime_j \geq m_j$ for all $j = 1, 2, \hdots, n$, we have 
$A_\text{XPLT}(\mathbf m^\prime) \geq A_\text{XPLT}(\mathbf m)$. In other 
words, it is never optimal to consider fewer candidates than what we have the
capacity for. This monotone property means we only have to consider extremal
points in the feasible set defined by the linear constraints above to solve the
integer program.


\subsection{A rich problem landscape}

Multi-stage screening is a rich and multi-faceted problem that can be studied
in various contexts. Above, we outlined a few dimensions of interest. To
summarize, we identified two \textbf{objectives}: multi-objective and
multi-fidelity optimization.  We must specify what aspects of the pipeline can
be optimized or controlled.  There are a variety of such pipeline
\textbf{design or decision variables}, but above we identified two: the
stage-wise allocations and, given a fixed allocation, which candidates to pass
between stages.  We outlined a few \textbf{policies} that each selects specific
settings of the design variables, including the XPLT policy that selects
allocations based on maximizing the expected reward.  We then identified some
\textbf{constraints}. First we distinguished between a closed-loop setting
versus a one-shot scenario, which is a constraint on the number of times we 
can execute the full real-world pipeline. We also discussed budget constraints to
limit the cost of one execution of the pipeline. 

By varying objectives, design variables, policies and constraints, we obtain a
rich problem landscape. In what follows, we shall focus on a small portion of
this landscape. Specifically, we turn our attention to selecting optimal
allocations which maximize the multi-fidelity reward function in the one-shot
setting, using the XPLT policy.  In section \ref{sec:results}, we detail some
simulation studies that measure the impact of various pipeline parameters on
the overall effectiveness of the screening procedure in this setting.

\section{A prior model for stage-wise scores}

As described above, we shall model the scores within a Bayesian context by
maintaining a multivariate probability distribution of the scores.  This
distribution captures predictions for the scores, an uncertainty quantification
for such predictions, and estimated correlation between different scores.
Starting from a prior distribution of such scores, we can integrate observed
scores for candidates selected by a policy using Bayes's Law. This results in
an updated posterior distribution on scores. Key to this procedure is the prior
distribution, which reflects our initial assumptions about the scores. Also
essential is the ability to efficiently sample scores from this prior
distribution, which we use to simulate the pipeline. Below, we detail these two
aspects.

We shall consider a flattened representation of the scores $\mathbf Y =
(Y_{i,j})$ so that $\mathbf Y$ is an $nm$-long vector whose $(jm + i)$-th entry
is $Y_{i,j} = f_j(x_i)$.  For our simulation studies, we shall assume that the
scores follow a  multivariate normal prior distribution, $\mathbf Y \sim
\mathcal N(0, \Gamma)$. That is, we assume that the scores on average are
centered at 0, which we view as a nominal baseline value.  Statistically
significant deviations away from this baseline value will represent ``good"
candidates. The covariance matrix $\Gamma$ is an $nm \times nm$ matrix, whose
entries describe the statistical relationship between the scores.

For our problem, we assume that the covariance is \emph{separable}. That
is, we can decompose $\Gamma$ into the Kronecker product $\Gamma = \Sigma
\otimes X,$ where $\Sigma$ and $X$ are $n \times n$ and $m
\times m$ covariances matrices, respectively. Here, $\Sigma$ represents
the relationship between stages and $X$ represents the relationship between
candidates.  Namely, the covariance $\text{Cov}(Y_{i,j}, Y_{k,\ell})$ is given 
by the product $\Sigma_{j, \ell}\cdot X_{i, k}$ of the $(j, \ell)$-th 
and $(i,  k)$-th entries of $\Sigma$ and $X$, respectively.  
That is, the relationship between the scores of different candidates at
different stages ``factors through" assumed inherent statistical
relationships between stages and candidates.

The assumption of separability is not a fundamental one, and we may imagine
problem settings where the covariance relationship between different scores is
non-separable. However, we choose this setting to perform our simulation
studies for several reasons. First, it is computationally efficient (see
Section \ref{sec:sampling}). The assumption also allows us to study the impact
of stages and candidates covariance structure independently.  Lastly, in
contrast to other methods that consider similar separable structure (such as
fitting a multi-task Gaussian Process \cite{bonilla2007multi} to model several
responses or tasks at once as a function of input variables), the assumption of
separability is more natural in our multi-fidelity setting.  Specifically, the
different tasks we are fitting a multivariate distribution to are different
models predicting the same response, albeit at different levels of fidelity. In
this setting, it is more natural to assume correlations between scores of two
different candidates persist to some degree between different fidelity levels,
i.e. the covariance is separable.

The optimal allocation problem is, in part, defined by this distribution.  It
serves as the prior distribution from which we sample ground truths to estimate
expected rewards. Different distributions will result in varying effectiveness
of the screening procedure. We wish to characterize the connection
between the distribution and the screening procedure's efficacy. To do this, we
generate distributions as follows.

Pick dimensions $d_x$ and $d_s$. Then independently and uniformly 
sample $m$ vectors $\left\{\mathbf x_1, ..., \mathbf x_m\right\}$ from the 
$d_x$-dimensional unit hypercube, $\mathbf x_i \sim \text{Unif} [0,1]^{d_x}$. 
With these samples, define the candidate covariance matrix $X$
\begin{equation}
X_{i,j} = \sigma_x^2 \exp\left[-\frac{\|\mathbf x_i - \mathbf x_j\|^2_2}{2\ell^2_x} \right], \quad i, j = 1, 2, \hdots,  m.
\label{eqn:prior_cand_cov}
\end{equation}
Here $\sigma_x, \ell_x$ are free parameters. Similarly define the stage
covariance $\Sigma$ matrix using $n$ independent samples 
$\left\{\mathbf s_1, \hdots \mathbf s_n\right\}$, where $\mathbf s_j \sim \text{Unif} [0,1]^{d_s}$:
\begin{equation}
\Sigma_{i,j} = \sigma_s^2 \exp\left[-\frac{\|\mathbf s_i - \mathbf s_j\|^2_2}{2\ell_s^2} \right], \quad i, j = 1, 2, \hdots,  n, 
\label{eqn:prior_stage_cov}
\end{equation}
for free parameters $\sigma_s, \ell_s$.  We call the samples $\mathbf X =
\left\{\mathbf x_i\right\}_{i=1}^m$ and $\mathbf S = \left\{\mathbf
s_j\right\}_{j=1}^n$ the \emph{latent representations} of candidates and
stages.  Different latent samples yield different covariance matrices, and
hence different prior distributions. The parameters $\sigma_x$ and $\sigma_s$
capture the magnitude of the scores assumed by the prior model.  Without loss
of generality, we shall assume both are equal to 1 for our simulation studies.

In reality, the covariance structure is problem-specific, yet the procedure
above provides a method for systematically varying the covariance structure in
a parametric way to study the impact of such structure on the effectiveness of
the screening pipeline.  Implicit in this generation procedure is the
assumption that the covariances are \emph{stationary} in the latent spaces.
That is, scores between different candidates or stages are considered
statistically similar if their latent representations are close together in the
latent space, and the notion of closeness is the same wherever the latent
representations are. This assumption of closeness in latent space as a proxy
for similarity between candidates is used in many generative design methods
(e.g. \cite{lim2018molecular}). Below, we study the impact of the latent
dimensions $d_s, d_x$  and length-scales $\ell_s, \ell_x$, which impact this
measure of closeness. There are other ways to sample or otherwise construct the
covariance matrices, such as sampling such matrices from an Inverse Wishart
distribution\cite{gelman1995bayesian}, or considering a non-stationary
covariance structure over the latent space.  Regardless of the sampling
procedure, we can perform a similar simulation analysis below given the
covariance matrices, $\Sigma, X$.

\subsection{Lazy stage-wise sampling}
\label{sec:sampling}
Once the prior $\mathcal N(0, \Gamma)$ is defined, we can use it to sample
ground truth for simulations.  Na\"ively, we may sample the entire set of $nm$
scores at once. This is wasteful because many candidates are screened out of
the pipeline, obviating the need for many of these scores in later stages. In
reality, only $M = \sum_{j=1}^n m_j$ scores are needed to simulate the
pipeline, and this sum is typically much smaller than $nm$.  Indeed, in many
applications, the stage costs increase exponentially, necessitating that
allocations decrease by an order of magnitude at each stage.  In this case, $M$
is essentially $km$, for some constant $k$.  

Computationally, na\"ively sampling from a fixed distribution $\mathcal N(0,
\Gamma)$ incurs a one-time $\mathcal O((nm)^3)$ computational cost and a
per-sample costs of 1) $nm$ samples from a univariate standard normal
distribution and 2) an $\mathcal O((nm)^2)$ matrix-vector multiplication. Such
per-sample costs are incurred for every simulation we run to obtain a reward
distribution for a specific allocation.  These costs are further multiplied by
the number of allocations considered when solving the integer program to
identify the optimal allocation. This is all to say that na\"ive sampling can
increase the costs of our simulation-based analysis drastically.  Improving how
we sample scores is therefore important for efficient simulation studies.

We can gain significant speed-ups in a straight-forward way using lazy
stage-wise sampling.  We sample scores for stage $j$ only when the
simulation of the pipeline arrives at this stage. At this point, we will know
the subset of candidates that had survived up to this stage, and
hence we shall only sample those scores. Suppose candidates $x_{i_1},
\hdots x_{i_{m_j}}$ have made it to stage $j$. Then we can sample the scores
$Y_{i_1, j}, \hdots, Y_{i_{m_j}, j}$ \emph{conditional} on all the previously
sampled scores from prior stages. By properties of the multivariate normal
distribution, these conditional scores are also multivariate normally
distributed. By additionally taking advantage of the separable structure of the
covariance matrix, we can perform the sampling in an efficient manner.  See
section \ref{sec:lazy_sample_algo} for details.

\subsection{Simulating a pipeline with an inter-stage policy}
With the ability to sample  ground truth scores, we can simulate the execution
of the pipeline under some pre-specified inter-stage policy and some fixed
allocation $\mathbf m$.  For the first stage, we evaluate all $m_1 = m$ initial
candidates $x_i$ by sampling scores $y_{i,1}$ for  $i = 1, ..., m$.  Using the
sampled scores, an inter-stage policy selects $m_2$ candidates to advance to
stage 2.  For example, the \emph{pure-exploitation inter-stage policy} simply
selects the $m_2$ candidates with highest stage-1 scores. 

Other inter-stage policies may consider the conditional distribution of the
\emph{final-stage} scores in light of the observed data. That is, conditioned
on the observations $\left\{y_{i,1}\right\}$, we can consider the distribution
of final scores $Y_{i, n}$ and select candidates for the second round based on
this distribution. This distribution is also multivariate normally distributed.
We may, for example use a delayed exploitation policy, where we select
candidates corresponding to the largest expected final scores under this
distribution. We can attempt to balance between exploration and exploitation as
well, using polices such as UCB.  We do not dwell on other inter-stage policies
in this paper. Instead, we will focus solely on the pure-exploitation
inter-stage policy.

Once the inter-stage policy selects the $m_2$ candidates for the second stage,
their scores for this stage are sampled, conditioned on the sampled scores
$\left\{y_{i,1}\right\}$ of the previous stage. Then the process repeats: based
on these scores, the appropriate number of candidates are selected for the next
stage per the pre-specified allocation. For these candidates, scores from the
next stage are sampled conditionally on the previously sampled scores, and so
on.  This is repeated until the evaluation of the last stage.  After simulating
the pipeline through the last stage, we compute the reward obtained for this
simulation. Then the procedure is repeated several times more, each simulation
using a different ground truth sample drawn from the prior.  These simulations
yield an empirical distribution of the reward $R(\mathbf m)$ for the given
allocation $\mathbf m$.  As described above, the meta-policy will use this
reward distribution to select an allocation $\mathbf m^\star$.

In practice, this decision $\mathbf m^\star$ would be then implemented, running
the real-life pipeline with this allocation. In the closed-loop design setting,
we would use the data obtained from this real-life execution of the pipeline to
update the Bayesian prior, which we would use in turn to select another
allocation to try. In the one-shot setting, we receive a single \emph{realized
reward} $r(\mathbf m^\star)$ from the sole execution of the pipeline.  Assuming
the prior is representative of the real ground truth, this realized reward
should be consistent with the distribution of rewards obtained via simulation.
That is, under the \emph{truth-from-prior} assumption (that the real ground
truth is sampled from the prior distribution) the distribution of realized
reward $r(\mathbf m^\star)$ is exactly that of $R(\mathbf m^\star)$, the
simulated reward distribution used to select $\mathbf m^\star$ in the first
place.  In the simulations below, we explicitly enforce this truth-from-prior
assumption.

\section{Results and discussion}
\label{sec:results}

In this section, we report results of simulation studies where we simulate
entire pipeline definition and optimization procedure: 
\begin{enumerate}
\item Sample latent representations of candidates and stages to define a 
prior distribution on scores.
\item Given the prior distribution, use an inter-stage policy to simulate the
pipeline via lazy stage-wise sampling and obtain samples of reward,
\item Solve the integer program defining the meta-policy to select an 
optimal allocation $\mathbf m^\star$, performing (2) above as an inner-loop.
\item  Characterize how well the meta-policy selects optimal allocations
through further out-of-sample simulations, here serving as a proxy for
real-world implementation of the selected allocation $\mathbf m^\star$,
\end{enumerate}
The outcome of step (4) above is a distribution for realized reward for the
allocation selected by the meta-policy.  Under the truth-from-prior
assumptions, the distribution of the realized reward is the same as the
distribution of the simulated rewards, and so we shall use the two
interchangeably. Specifically, below we will simply refer to the expected
optimal rewards, 
\[ \begin{aligned}
\overline R(\mathbf m^\star) 
    &= \mathbb E[R(\mathbf m^\star)] \\
    &= \mathbb E[r(\mathbf m^\star)], \text{ (under truth-from-prior assumption.)}
\end{aligned}\]

In the simulation studies below, we characterize how different pipeline
parameters impact the screening procedure's efficacy, as measured by the
expected optimal rewards obtained.  We compare against a baseline policy, which
we refer to as the \emph{random} policy. Under the random policy, we expend as
much of the budget $C_\text{max}$ by testing candidates solely on the final
stage, without any screening. We select candidates uniformly at random, with
replacement. Once the budget is spent on these random final-stage trials, we
calculate the reward in the same manner as the screening procedure: select the
maximum score from among the randomly selected candidates. Repeating this
exhaustive random selection several times results in an expected reward, $\overline R$,
for the random policy.

\begin{table}
\centering
\begin{tabular}{c|l|c} \hline
Parameter & Description & Value \\ \hline
$m$ & Number of initial candidates & 500 \\
$C_\text{max}$ &Total budget & 2500 \\
$d_x$ & Dimensionality of latent candidate representation & 8 \\
$\ell_x, \sigma_x$ & Parameters to candidate covariance function 
    (Eqn. \eqref{eqn:prior_cand_cov}) & 1, 1\\
$d_s$ & Dimensionality of latent stage representation & 1 \\
$\ell_s, \sigma_s$ & Parameters to stage covariance function 
    (Eqn. \eqref{eqn:prior_stage_cov}) & 0.2, 1\\ \hline
\end{tabular}
\caption{Base parameter values used throughout the simulation studies.}
\label{tbl:base_params}
\end{table}

We consider a three-stage pipeline, with the three stages representing a cheap
machine-learning material property predictor, an expensive physics-based model,
and a physical experiment, respectively. We assign  nominal exponentially
increasing costs $\mathbf c = (1, 10, 100)$ for the stages. Unless otherwise
stated, the simulations use a set of base parameters listed in Table
\ref{tbl:base_params}.

\subsection{Effect of budget $C_\text{max}$ and number of initial candidates $m$}
\label{sec:simstudies}
One reason to carry out simulation studies is to assess the sensitivity of the
screening procedure on specific parameters. The supplemental information
(Section \ref{sec:si_simstudies}) provide results for several such simulation
results, in which specific parameters were varied. Here we focus on two for
illustration.

\begin{figure}
\centering
\begin{subfigure}[b]{0.35\textwidth}
\centering
\includegraphics[width=\textwidth]{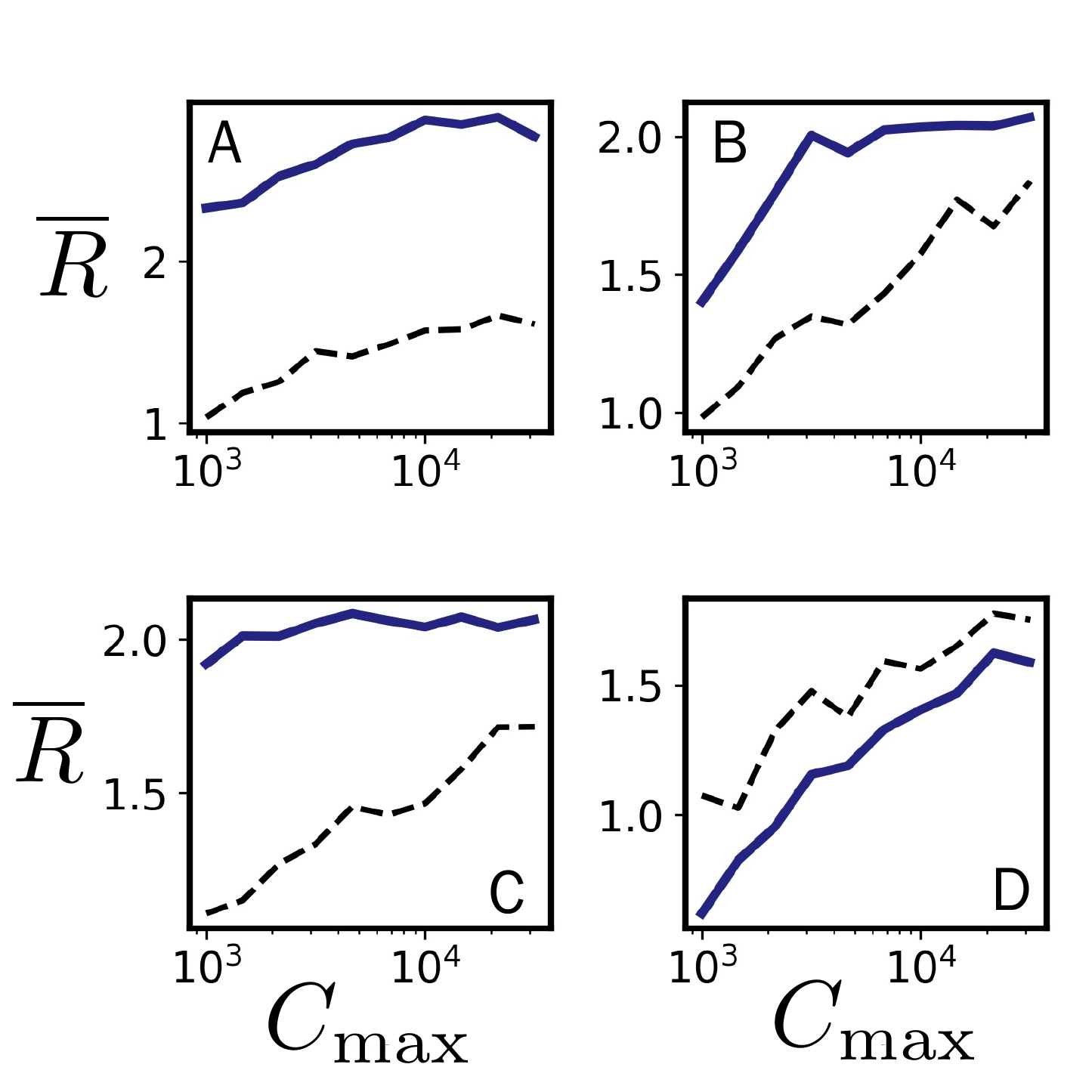}
\caption{Expected rewards vs. $C_\text{max}$}
\label{fig:c_max}
\end{subfigure}
~
\begin{subfigure}[b]{0.35\textwidth}
\centering
\includegraphics[width=\textwidth]{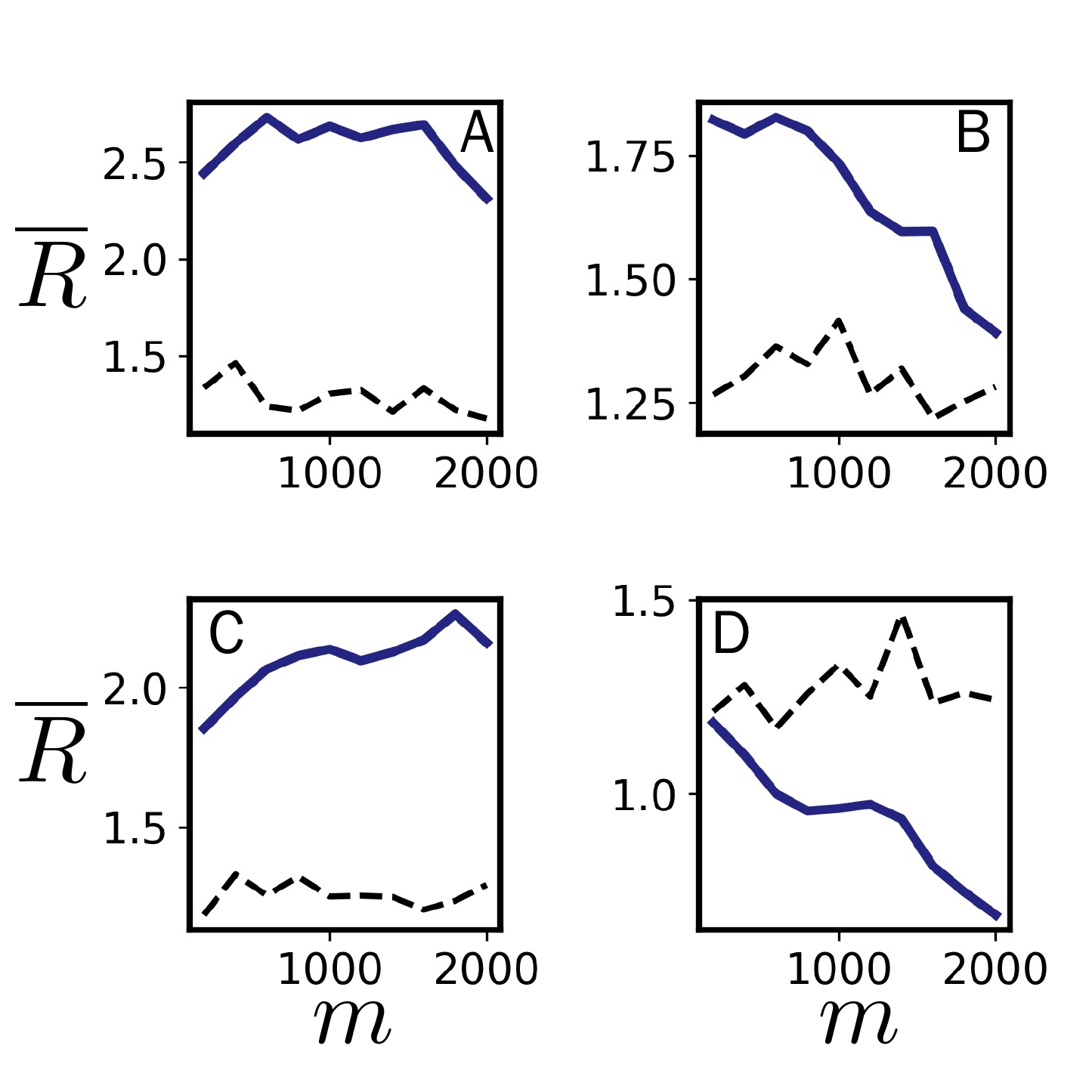}
\caption{Expected rewards vs. $m$}
\label{fig:m}
\end{subfigure}

\caption{Expected rewards vs. specific pipeline parameters. Solid blue lines
are the expected optimal rewards obtained from the XPLT-based screening
procedure, while the dashed black line show the expected rewards obtained under
the random (non-screening) policy. Each sub-panel (A, B, C, and D) corresponds
to a different prior distribution, and hence a different problem.
Specifically, the prior distributions differ with respect to stage-wise
similarities, as described in Section \ref{sec:effect_prior_dist}.}
\end{figure}

Figure \ref{fig:c_max} shows the expected reward obtained as a function of the
budget $C_\text{max}$, which is varied from 1000 to 50,000. We note that a
budget of 50,000 is exactly the cost to exhaustively evaluate the $m = 500$
candidates once during the final stage, each evaluation incurring a cost $c_3 =
100$. The solid blue lines plot the expected optimal reward under the XPLT
meta-policy, while the dashed black line shows the expected rewards under the
random policy.  The four plots each correspond to four different choices of
priors, and hence represent four different problems, each with unique
statistical relations between candidates and stages. We note that in general,
the performance of both policies increase as the computational budget
increases. Yet, depending on the distribution, the marginal returns obtained
from increasing the budget varies. For example, in panel C, there is little
marginal improvement to the expected reward even as we increase the budget by
an order of magnitude. In panel B,  we obtain good marginal improvement only up
to $C_\text{max} = 3000$.  The bottom-right panel (panel D) shows it is
possible for the screening procedure to consistently perform worse than random
for some problem settings.

Figure \ref{fig:m} shows the expected reward obtained versus the number of
initial candidates $m$, varying this number between 500 and 2000. We note that
the performance of the random policy is not sensitive to this value. With the fixed
budget $C_\text{max} = 2500$ and a final-stage cost $c_3 = 100$, the random
policy can try 25 out of $m$ total candidates. As the simulations suggest, the
chance that one of 25 selected candidates yields a large reward is relatively
stable (and low) for $m$ between 500 and 2000. In contrast, the expected
optimal rewards obtained from the XPLT policy shows varying behavior depending
on the prior distribution. In panel A, the expected optimal rewards is
relatively insensitive to $m$, while in panel C, these rewards increase with
$m$. Panels B and D show the opposite behavior, where the performance decreases
with respect to $m$.  As with the previous study, we see that the
screening procedure performs worse than random in panel D. 

It is interesting to note the different trends in expected optimal reward
versus different problem parameters or distributions of the ground truth, and
how such trends compare to random. The above simulations and those in
\ref{sec:si_simstudies} make clear that the analysis is not straight forward
and that there is no one general ``rule-of-thumb" when it comes to setting these
parameters.  Therefore, in practice, it is important to assess the problem
setting (in the above: which of A, B, C, or D) prior to executing a screening
pipeline or even selecting policy parameters using \emph{ad hoc} or heuristic
methods. Simulations offer a systematic way of studying such effects.

\subsection{Effect of prior distribution}
\label{sec:effect_prior_dist}

\begin{figure}
\centering
\includegraphics[width=0.5\textwidth]{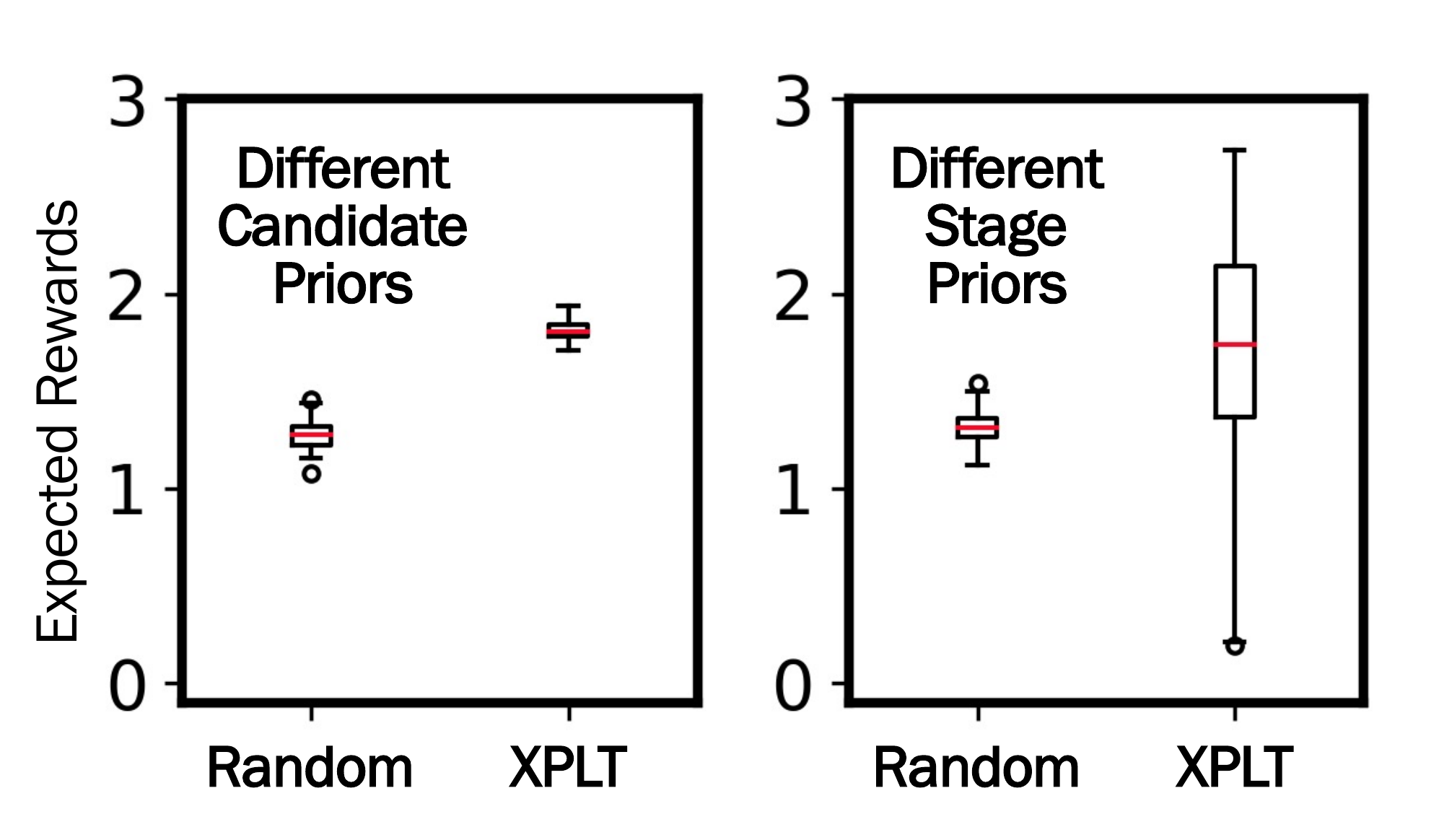}
\caption{Expected reward distributions over different priors.}
\label{fig:prior_sensitivity}
\end{figure}

As we see above, the trends of expected rewards could depend on the prior
distribution, since the real ground truth is sampled from this distribution in
these truth-from-prior studies. It is therefore desirable to characterize
this dependence further, and see if there are any features of a distribution 
or problem that predict effectiveness of the screening procedure.  The prior
distribution is characterized by the latent samples of candidates, $\mathbf X =
\left\{\mathbf x_i\right\}_{i=1}^m$, of stages $\mathbf S =
\left\{\mathbf s_j\right\}_{j=1}^n$, and parameters such at the dimensionality
of such samples $d_x$ and $d_s$, and the covariance function hyperparameters
$\ell_x, \ell_s, \sigma_x$ and $\sigma_s$.  The SI details simulation studies
for these parameters. Here, we study the impact of the samples $\mathbf X$ and 
$\mathbf S$. 

By considering several such latent samples, we obtain distributions of expected
rewards, which are shown in Figure \ref{fig:prior_sensitivity}.  The left panel
shows the induced distribution on expected reward as we vary the latent
samples on candidates $\mathbf X$, for both the XPLT policy and the random
baseline policy. As we see, the distributions for both policies are insensitive
to the specific choice of latent samples, typified by the tight induced
distributions.  In contrast, the right panel shows these distributions as we
vary stage latent samples $\mathbf S$. Here, we see that the expected
reward is highly sensitive to the actual latent samples of stages, implying
that stage-wise covariance structure significantly impacts the overall
effectiveness of the screening procedure. We even see that a significant
proportion of sampled priors result in the screening procedure performing worse
than random.

\begin{figure}
\centering
\includegraphics[width=0.7\textwidth]{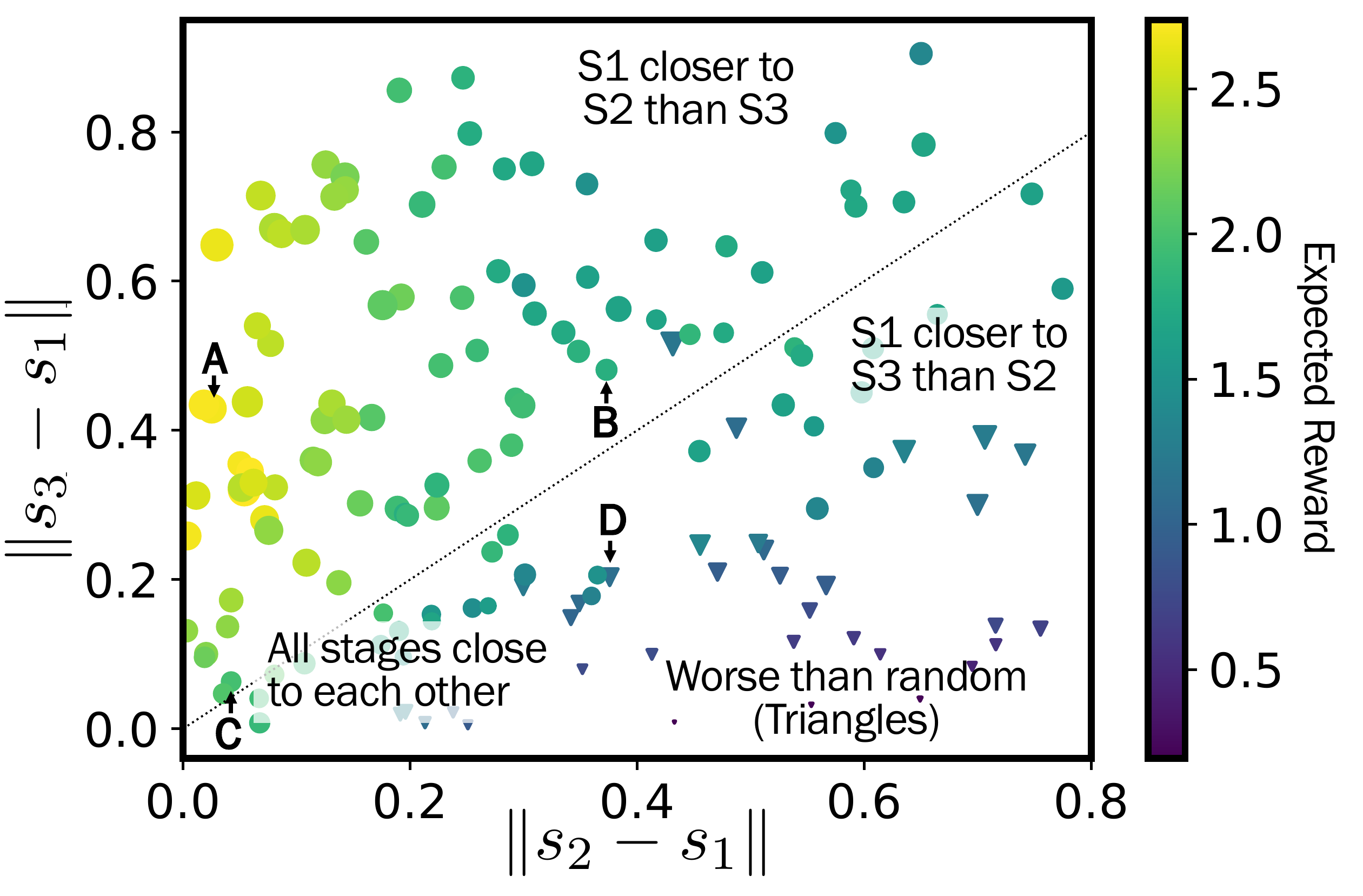}
\caption{Similarities between stages dictate expected rewards.  This plot shows
expected realized rewards for different prior distributions.  Locations of the
plotted points indicate the distance between first and second stages, and first
and third stages.  Color indicates the expected realized rewards, while size is
inversely proportional to the variance of the realized rewards.  Inverted
triangles indicate cases where the expected reward obtained via screening was
less than that obtained via random sampling.}
\label{fig:heatmap}
\end{figure}

To further explore this phenomenon, we consider the optimal rewards obtained by
the XPLT policy for several samples of $\mathbf S = (\mathbf s_1, \mathbf s_2,
\mathbf s_3)$, and hence for several different priors distributions.  We
characterize these rewards with respect to two geometric features derived from
the latent stage representation: the distance between the latent
representations of Stage 1 and Stage 2 $\|\mathbf s_2 - \mathbf s_1\|$
and the distance between Stage 1 and Stage 3 $\|\mathbf s_3 - \mathbf
s_1\|$. These distances typify the statistical relationship between the stages.
This data is displayed in Figure \ref{fig:heatmap}.

Each dot in the figure corresponds to a specific sample of latent stages i.e. a
specific distribution for ground truth. Dots are located at coordinates
$(\|\mathbf s_2-\mathbf s_1\|, \|\mathbf s_3 - \mathbf s_1\|)$. For example,
points in the lower-left corner correspond to the case where the three stages
are statistically similar, while dots in the top-left corner correspond to the
case where Stages 1 and 2 are  similar, but the Stages 1 and 3 are not. The
color of a dot indicates the expected optimal reward given the specific prior
distribution corresponding to that dot. The size of the dot is inversely
proportional to the variance of the reward distribution. Smaller dots are cases
where there is a wide distribution of rewards associated to the selected
optimal allocation.  Circles indicate the cases where the XPLT policy performs
better than random given the same prior distribution, while inverted triangles
show the cases where it performs worse than random. The points labeled A, B, C
and D correspond to the specific prior distributions used in the
single-parameter simulation studies in Section \ref{sec:simstudies}
above and in the SI.

The dashed line indicates the points where Stage 3 and Stage 2 are both
equidistant to Stage 1.  Points above this line refer to cases where Stage 1 is
more statistically similar to Stage 2 than to Stage 3.  Points below this line
correspond to cases where Stage 1 is more similar to Stage 3 than to Stage 2.
In this region, the statistical similarities between stages do not respect the
order in which the stages are executed in the pipeline.  From the plot, the
screening procedure is not effective in this region, often performing worse
than random.  There is a general improvement as Stage 2 and Stage 1 become more
similar, i.e.  as we move right to left along the horizontal axis.
Interestingly, this trend does not hold with respect to Stage 3. From the plot,
there appears to be an optimal \emph{dissimilarity} between Stages 1 and 3 that
yield maximal expected rewards. For example, among the samples, Point A
corresponds to the largest expected reward. For this point the distance between
$\mathbf s_3$ and $\mathbf s_1$ is around 0.42. 

This optimal dissimilarity for this model system again indicates the need to
perform such simulations for real systems in practice. In this specific case,
tuning similarity between stages to match the optimal settings identified
through such simulation studies could mean tuning the accuracy of a ML model,
increasing or decreasing simulation run times or resolutions, or introducing
randomness (i.e. exploration) in selecting candidates to advance between
stages. Having a proper assessment of the impact of such pipeline parameters
could then result in more efficient and effective screening.

\subsection{Throughput and optimal reward}
\begin{figure}
\centering
\includegraphics[width=0.5\textwidth]{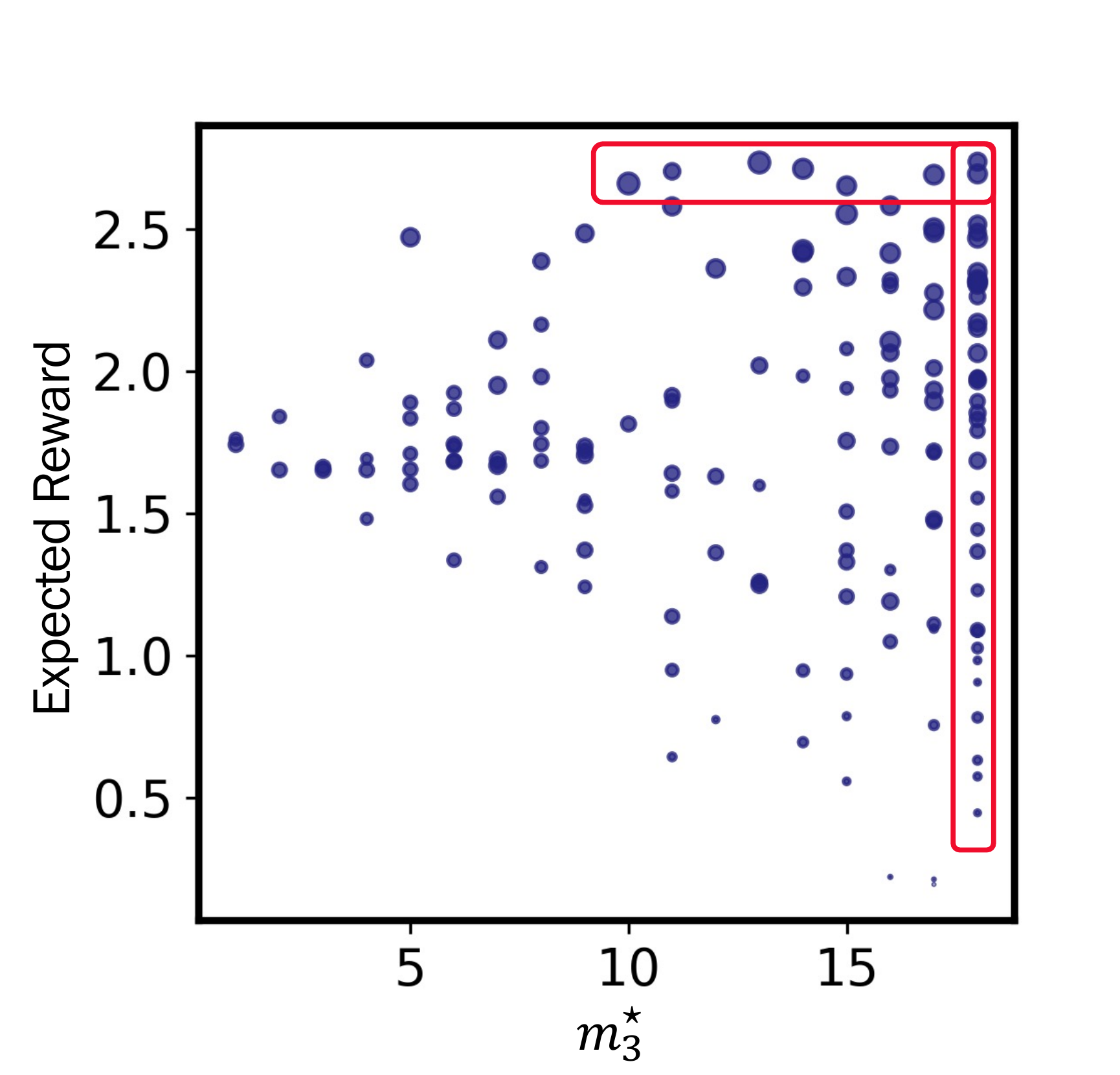}
\caption{Throughput correlates with expected optimal reward. This figure plots
expected optimal reward versus final-stage allocation $m^\star_3$ of the
optimal allocation $\mathbf m^\star$ selected by the XPLT policy. Each dot
corresponds to a different prior distribution. The size of each dot is
inversely proportional to the variance of the reward distribution $r(\mathbf
m^\star)$.} \label{fig:throughput}
\end{figure}

Another pipeline study \cite{woo2021optimal} examined the pipeline problem
under different design variables and objectives. There, authors considered
score thresholds (rather than explicit allocations) as the design variables.
They also considered the objective of optimizing pipeline \emph{throughput} --
the number of candidates that satisfy all thresholds to make it past the final
stage, under similar cost constraints as those detailed above. In this setting,
the pipeline is used to maximize the number of good candidates passing through
the screening procedure.  While not the same problem as considered here, we are
nevertheless able to draw some connections between our work and this optimal
throughput problem.

Figure \ref{fig:throughput} plots the expected optimal reward
$\overline{R}(\mathbf m^\star)$ versus the final stage allocation $m^\star_3$,
for different sampled prior distributions.  Each dot corresponds to the
expected optimal reward and final stage allocation obtained for a specific
prior distribution. As before, the size of the dot is inversely proportional to
the variance of the reward distribution. The red rectangles call attention to
settings with high expected rewards (top), and those that have large final
allocation (right). The final allocation is the measure for throughput. While
in our study we optimized for expected rewards rather than throughput, we
observe from the plot that the two quantities are correlated

Namely, we note that for the settings resulting in the highest expected reward
(top box), the final stage allocation $m^\star_3$ was around 10 or larger.
Given the final stage cost $c_3 = 100$ and budget $C_\text{max} = 2500$, this
means that optimal allocations yielding highest expected rewards were expending
around 40 to 60\% of the budget on final stage trials. We note that for cases
where the optimal final allocation is large (right box), there is a spread of
both high and low expected rewards. Yet, the concentration is biased toward the
high expected reward regime. Moreover the high-throughput / low expected reward
cases are associated with large variance of reward distribution. These results
show that while the relationship between throughput and expected rewards may be
complex, there are broad agreements between the two optimization objectives,
though there are cases where this agreement may not hold.

\section{Conclusion}

In this paper, we defined a model and algorithm for simulating pipelines,
and identified places where such methods can be used to perform decision-making
under uncertainty via simulation-based policies. We codified how various
different pipeline objectives fit inside the common framework of optimizing
expected rewards under uncertainty, and alluded to a multi-faceted, rich
problem set that belies the ubiquitous experimental protocol of multi-stage
screening.

Examining one part of this problem landscape, we focused on analyzing a popular
mode for screening: simply selecting the top $m_j$ candidates from the previous
stage to pass onto the next stage. We showed how to optimize expected rewards
with respect to a specific design parameter choice, allocation, in a one-shot
setting. In this specific context, our simulation studies showed a variety of
trends as we attempted to characterize the effectiveness of the screening
procedure with respect to pipeline or problem parameters. We most
importantly described  how stage-wise covariance structure played a significant
role in determining this efficacy. We demonstrated how simulated studies
allowed us to \emph{quantify} the performance, its dependence on parameters,
and in comparison with a random baseline policy. 

We hope to have demonstrated the need for proper modeling and analysis
in multi-stage screening  pipeline optimization, decision-making and parameter
calibration.  We believe that the models, policies and simulation methods
presented in this paper provide important tools in this analysis. 

\section{Acknowledgments}
This material is based upon work supported in part by the National Science
Foundation under Grant No. 1950796, NSF REU Site:\emph{``Data-driven Materials
Design."} The work was also supported in part by the Brookhaven National
Laboratory Directed Research and Development (LDRD) Grant No. 21-044.  We thank
Bill Bauer and Erik Einarsson for organizing the REU Site. We thank Byung-Jun
Yoon, Nathan Urban and Frank Alexander for helpful discussions.

\section{Conflict of interest statement}
On behalf of all authors, the corresponding author states that there is no
conflict of interest.

\bibliographystyle{ieeetr}
\bibliography{mybib}

\begin{appendix}

\renewcommand{\thesection}{SI.\arabic{section}}

\section{Iterative algorithm for lazy stage-wise sampling}
\label{sec:lazy_sample_algo}

In this section, we introduce the equations to express the parameters of the
relevant marginal and conditional normal distributions used to perform lazy
stage-wise sampling of scores, under separable covariance assumption.  As
before, we shall let $n$ be the number of stages, and $m$ be the number of
initial candidates. Let 
\[\mathbf Y = \begin{pmatrix} \mathbf Y_1 \\  \mathbf Y_2 \\ \vdots \\ \mathbf Y_n\end{pmatrix} \]
be the $nm$-dimensional vector of scores, expressed in block form so that
$\mathbf Y_j \in \mathbb R^m$ represents the $m$ scores for all candidates at
stage $j$. Concatenating all but the first block, we may write
\[ \mathbf Y = \begin{pmatrix} \mathbf Y_1 \\ \mathbf Y_{2:} \end{pmatrix},\]
with the concatenation is denoted as $\mathbf Y_{2:} = (\mathbf Y_2, \hdots, \mathbf Y_n)^T$.

Suppose $\mathbf Y \sim \mathcal N(\boldsymbol \mu, \Gamma)$ is multivariate
normal distributed, with the mean vector $\boldsymbol \mu \in \mathbb R^{nm}$ 
and covariance matrix $\Gamma \in \mathbb R^{nm \times nm}$ given in block
form as:
\[ \boldsymbol \mu = \begin{pmatrix} \boldsymbol \mu_1 \\ \boldsymbol \mu_{2:} \end{pmatrix}, \quad\quad
   \Gamma  = \begin{pmatrix} \Gamma_{1,1}&\Gamma_{1,2:} \\ \Gamma_{2:,1} & \Gamma_{2:,2:} \end{pmatrix}. \]
Here, $\Gamma_{1,1}$ is the $m \times m$ covariance matrix between all
candidates for the first stage, while $\Gamma_{1,2:} = \Gamma_{2:,1}^T$ is the
$m \times (nm-m)$ matrix of cross-covariances between the scores of all
candidates at the first stage, and at later stages. $\Gamma_{2:,2:}$ is the
$(nm-m)\times(nm-m)$ covariance matrix for the scores of all candidates at all
other stages except for the first.  As before, we view $\Gamma$ as separable 
and write $\Gamma = \Sigma \otimes X$. We can similarly write the stage
covariance matrix $\Sigma$ in block-form as:
\[ \Sigma = \begin{pmatrix} \Sigma_{1,1}&\Sigma_{1,2:} \\ 
            \Sigma_{2:,1} & \Sigma_{2:,2:} \end{pmatrix}, \]
so that $\Sigma_{1,1}$ is a scalar, $\Sigma_{2:,1} = \Sigma_{1,2:}^T$ is
an $(n-1)\times 1$ vector and $\Sigma_{2:, 2:}$ is the $(n-1)\times(n-1)$
sub-covariance matrix for stages 2 through $n$. It follows that $\Gamma_{1,1} =
\Sigma_{1,1} \otimes X$, $\Gamma_{2:,1} = \Sigma_{2:,1}\otimes X$ and
$\Gamma_{2:, 2:}= \Sigma_{2:, 2:} \otimes X$ 

We may obtain a sample of just the first stage scores $\mathbf Y_1$ by noticing
the marginal distribution is also normal. That is $\mathbf Y_1 \sim \mathcal
(\boldsymbol \mu_1, \Gamma_{1,1})$. Suppose $X$ admits a Cholesky factorization
$X = LL^T$. Then $\Gamma_{1,1}$ admits a Cholesky factorization 
$\Gamma_{1,1} = (\sqrt{\Sigma_{1,1}}L)(\sqrt{\Sigma_{1,1}}L)^T$. We can use this 
to sample first-stage scores:
\[\mathbf{y}_1 = \boldsymbol \mu_1 + \sqrt{\Sigma_{1,1}} L \mathbf z,\]
where $\mathbf z \sim \mathcal N(0, I_m)$ is a sample of $m$ independent, 
standard normal random variables.

Once the first stage scores are sampled, we can update the distribution of 
scores for later stages, conditional on the sampled scores. This conditional
distribution is also normal, namely 
$\mathbf Y_{2:} | \mathbf Y_1 = \mathbf y_1 \sim \mathcal N\left(\boldsymbol \mu_{2:}^\prime, \Gamma^\prime_{2:,2:}\right),$
where
\[
\begin{aligned}
    \boldsymbol \mu_{2:}^\prime &= \boldsymbol \mu_{2:} + \Gamma_{2:,1} \Gamma_{1,1}^{-1}(\mathbf y_1 - \boldsymbol \mu_1), \\ 
        &= \boldsymbol \mu_{2:} + (\Sigma_{2:,1} \otimes L) (\Sigma_{1,1} \otimes L)^{-1}(\mathbf y_1 - \boldsymbol \mu_1), \\
        &= \boldsymbol \mu_{2:} + \frac{1}{\Sigma_{1,1}} \left(\Sigma_{2:,1}\otimes (\mathbf y_1 - \boldsymbol \mu_1)\right), \\
    \Gamma_{2:,2:}^\prime &= \Gamma_{2:,2:} - \Gamma_{2:,1}\Gamma_{1,1}^{-1}\Gamma_{1,2:}. \\
                     &= \left(\Sigma_{2:,2:} - \Sigma_{2:,1}\Sigma_{1,1}^{-1}\Sigma_{1,2:}\right) \otimes X, \\
                     &= \left(\Sigma_{2:,2:} - \frac{\Sigma_{2:,1}\Sigma_{2:,1}^T}{\Sigma_{1,1}}\right) \otimes X, \\
                     &= \tilde \Sigma \otimes X,
\end{aligned}
\]
Here, $\tilde \Sigma = \Sigma_{2:,2:} - \Sigma_{2:, 1}\Sigma_{2:,1}^T / \Sigma_{1,1}$

In addition to sampling scores for the first stage, and updated the
distribution for later stages, the inter-stage policy will choose some
candidates that will not pass through to the second stage, and hence we do not
need to sample scores for these entries at later stages. Namely, we shall
filter out $\Delta m = m_2 - m_1$ candidates.  Let $\tilde L$ be the $m_2
\times m$ matrix resulting from removing the appropriate $\Delta m$ rows of $L$
and let $\tilde{\boldsymbol \mu}$ be the $m_2\cdot(n-1)$-length vector resulting
from removing the appropriate $(n-1)\cdot\Delta m$ rows from $\boldsymbol \mu_{2:}$
corresponding to the filtered candidates for each of the remaining $n-1$ stages.

At the end of this procedure, we have the updated conditional distribution: 
\[
\begin{aligned}
    \tilde{\mathbf Y} &\sim \mathcal N(\tilde{\boldsymbol \mu}, \tilde \Gamma),\\
    \tilde \Gamma &= \tilde \Sigma \otimes \Tilde X,\\
    \tilde X &= \tilde L \tilde L^T.
\end{aligned}
\]
which describes the score distribution for the remaining $m_2$ candidates, 
for the remaining $n-1$ stages. The procedure above can be repeated to sample
scores for the second stage and update the distribution for the third, and so
on. As we see above, the sampling and conditioning operations depend on our 
ability to maintain and update the mean score vector $\boldsymbol \mu$, the
stage covariance matrix $\Sigma$ and the Cholesky factor of the candidate
covariance matrix $L$ as we simulate the pipeline. At stage $j$ of the
pipeline, these updated matrices are of size $m_j\cdot n$, $(n-j+1) \times
(n-j+1)$ and $m_j\times m$, respectively, a dramatic departure from the
$(nm)\times(nm)$ scaling required for na\"ive sampling.

\section{Further simulation studies}
\label{sec:si_simstudies}

\begin{figure}[h]
\centering

\begin{subfigure}[b]{0.4\textwidth}
\centering
\includegraphics[width=\textwidth]{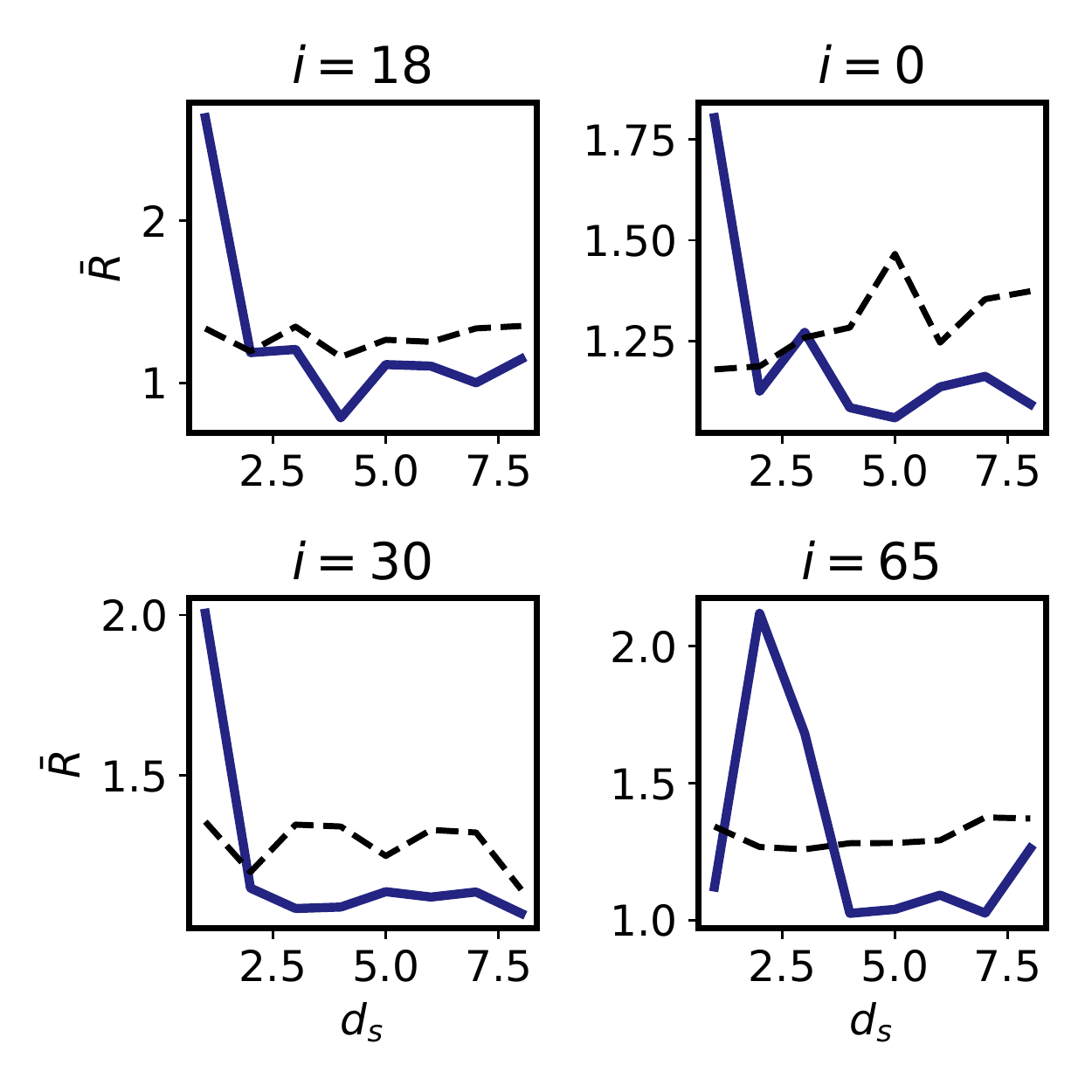}
\caption{$d_s$}
\end{subfigure}
~
\begin{subfigure}[b]{0.4\textwidth}
\centering
\includegraphics[width=\textwidth]{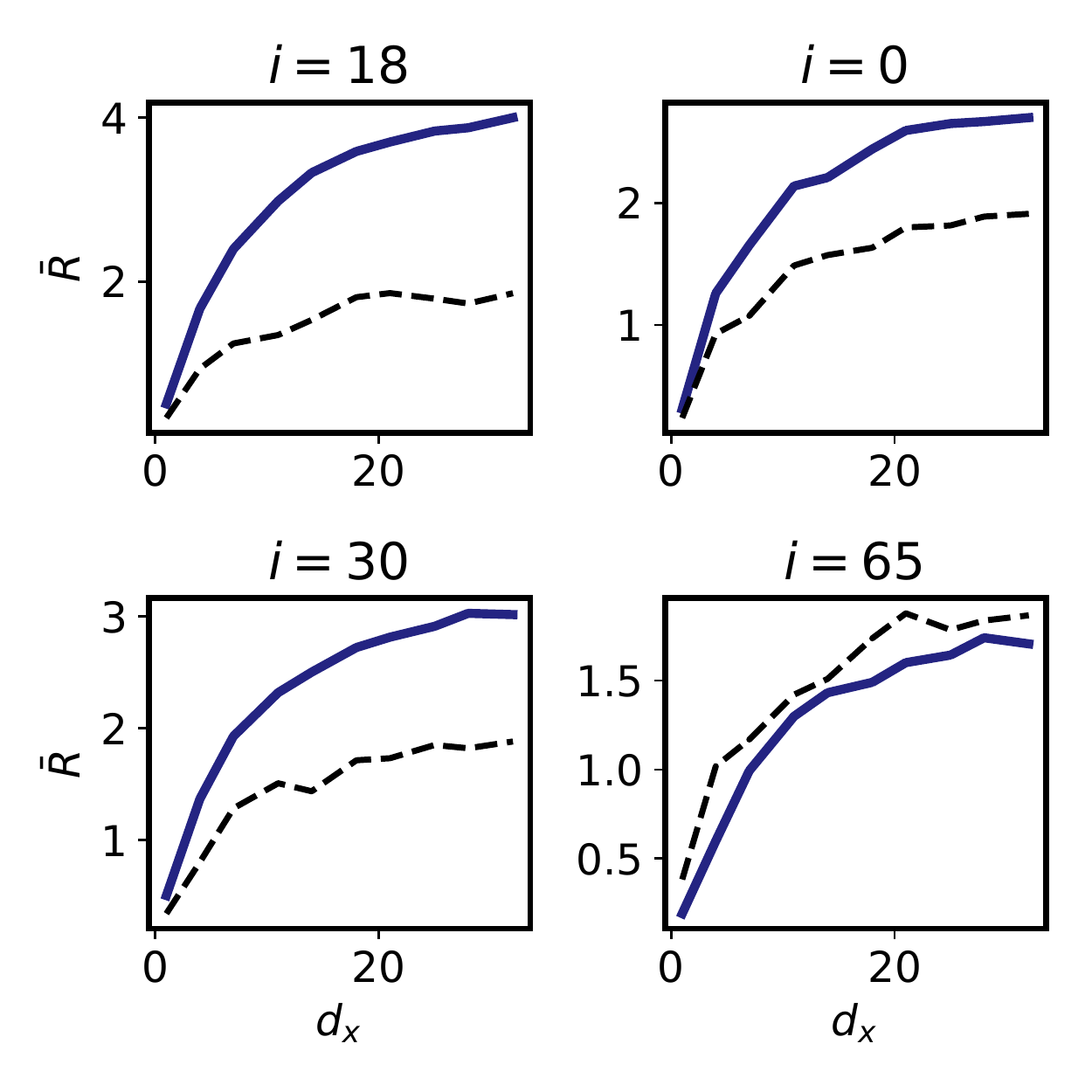}
\caption{$d_x$}
\end{subfigure}

\begin{subfigure}[b]{0.4\textwidth}
\centering
\includegraphics[width=\textwidth]{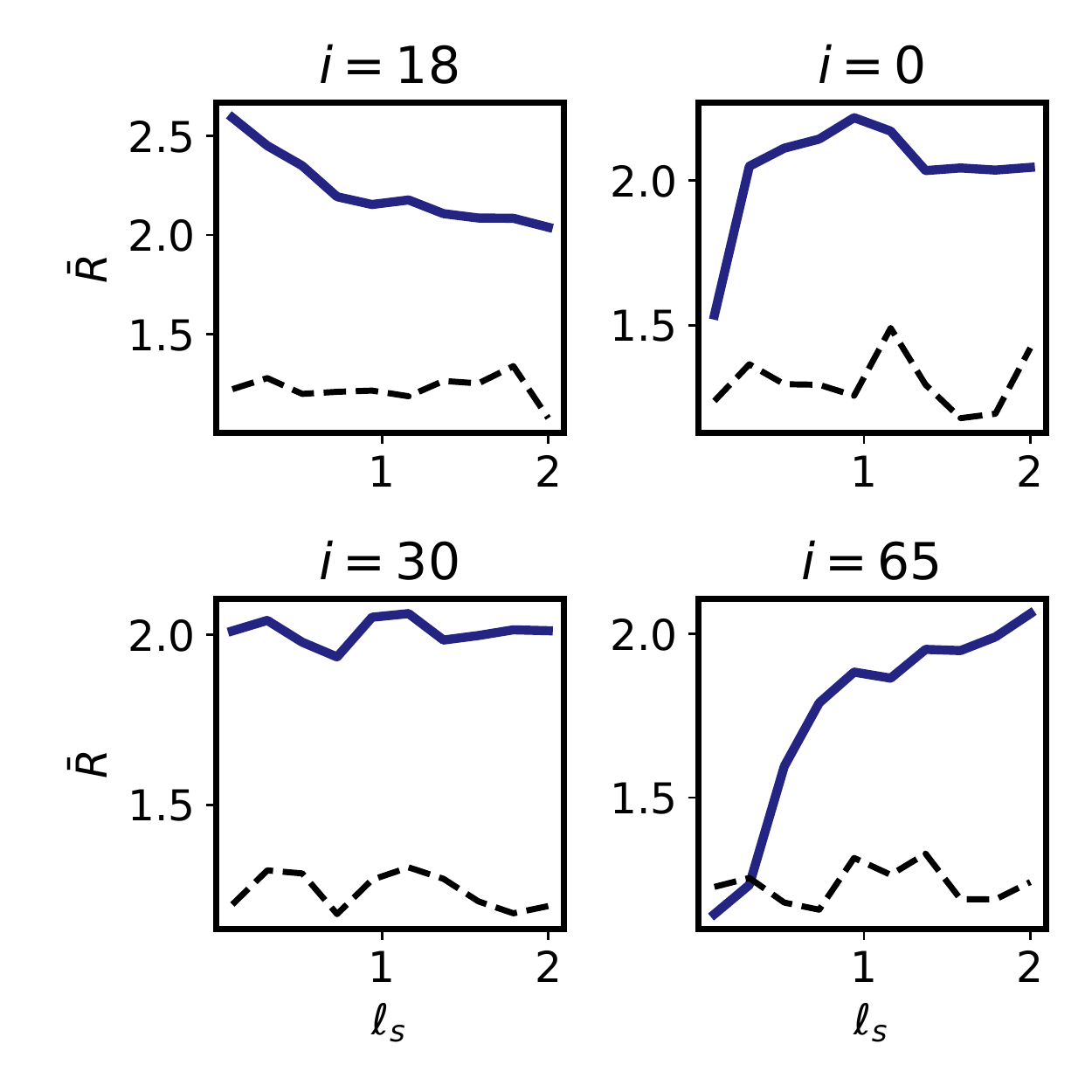}
\caption{$\ell_s$}
\end{subfigure}
~
\begin{subfigure}[b]{0.4\textwidth}
\centering
\includegraphics[width=\textwidth]{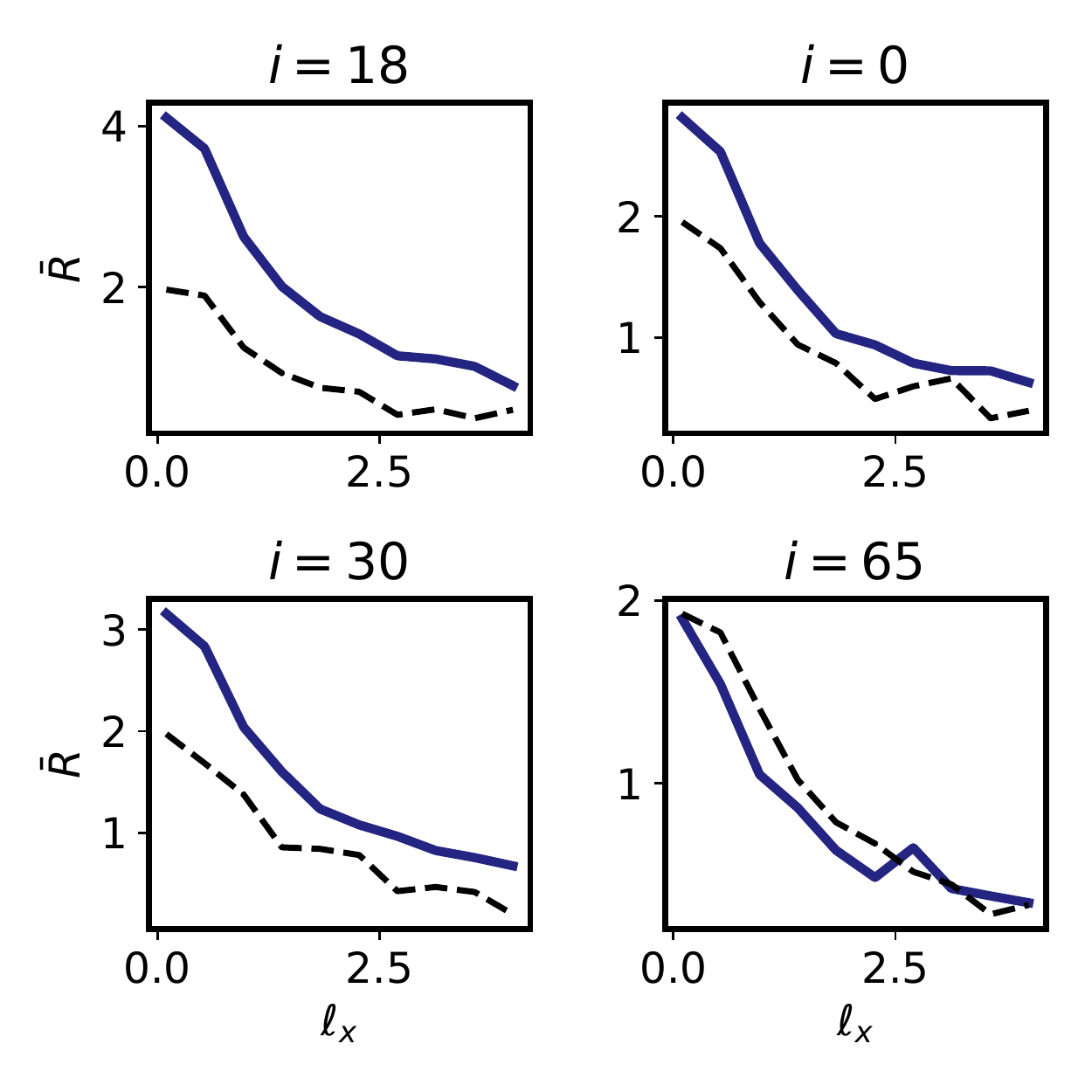}
\caption{$\ell_x$}
\end{subfigure}

\caption{Single parameter simulation studies for covariance prior parameters.
The order of the panels correspond to points A, B, C, and D as in the main
text.}

\end{figure}

\begin{figure}[h]
\centering
\includegraphics[width=0.7\textwidth]{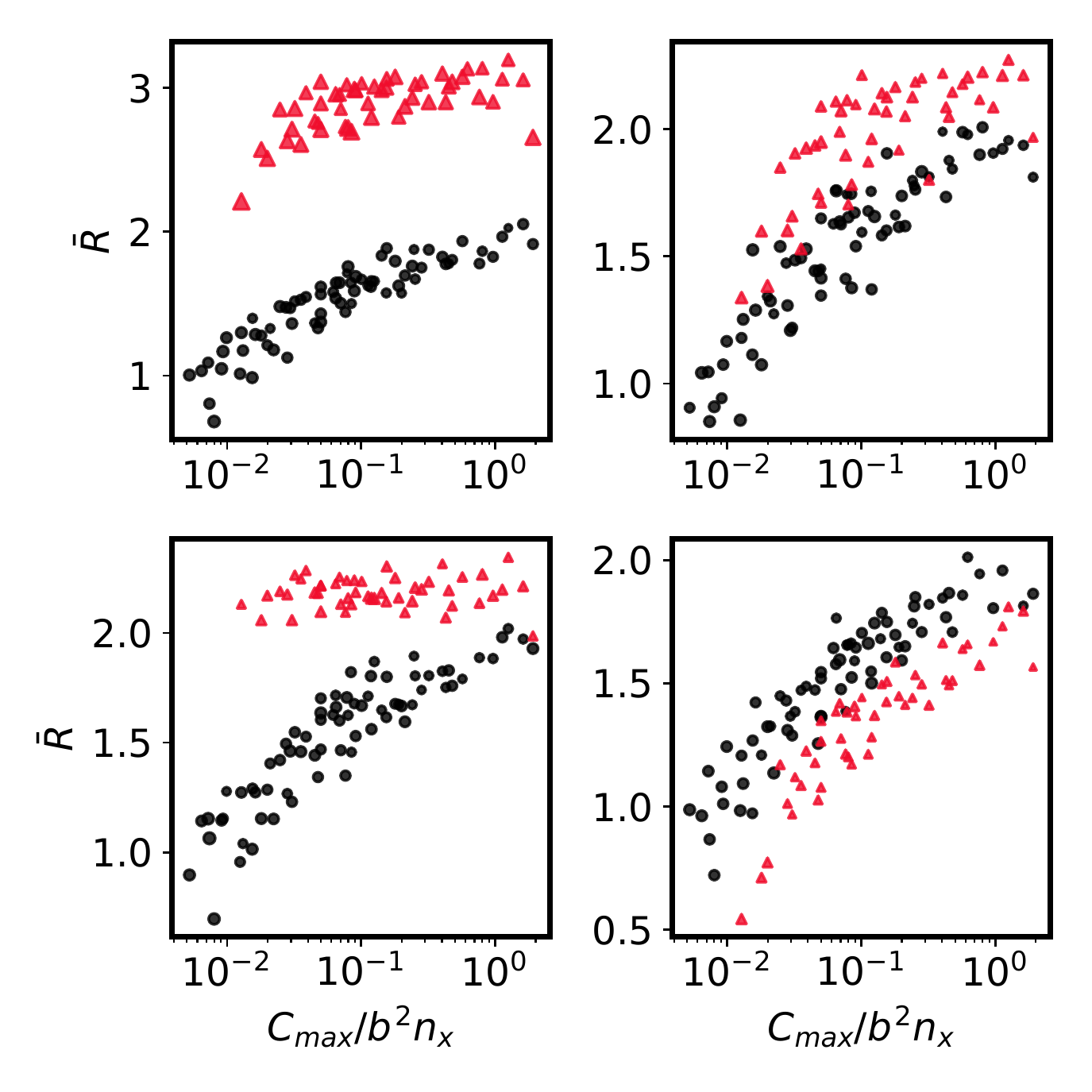}
\caption{ We parameterize cost as $\mathbf c = (b^0, b^1, b^2)$. The quantity
$C_\text{max}/b^2m$ is the proportion between the given budget $C_\text{max}$
and the cost to exhaustively try all candidates $b^2m$. We plot this versus
expected rewards for various choices of $C_\text{max}, b$ and $m$.  }
\end{figure}

\end{appendix}

\end{document}